\documentclass[11pt,reqno]{article}
\usepackage{amssymb, amsmath, amsthm, amsfonts, amscd}
\usepackage[english]{babel}
\usepackage{graphicx}
\newtheorem{theorem}{Theorem}[section]
\newtheorem{cor}[theorem]{Corollary}
\newtheorem{lemma}[theorem]{Lemma}

\newcommand{\nm}{\noalign{\smallskip}}

\def\R{\mathbb{R}}

\def\ep{\varepsilon}

\def\ue{u_{\varepsilon}}
\def\Ce{C_{\varepsilon}}
\def\fe{\varphi_{\varepsilon}}
\def\fo{\varphi_1}
\def\ft{\varphi_2}

\def\fte{\widetilde{\varphi}_\ep}

\def\XXint#1#2#3{{\setbox0=\hbox{$#1{#2#3}{\int}$}
     \vcenter{\hbox{$#2#3$}}\kern-.5\wd0}}

\newcommand{\Om}{\Omega}
\newcommand{\bX}{\mathbf{X}}

\newcommand{\vp}{\varphi}

\newcommand{\p}{\partial}
\newcommand{\pd}[2]{\frac {\p #1}{\p #2}}
\newcommand{\ds}{\displaystyle}
\newcommand{\eqnref}[1]{(\ref {#1})}

\newcommand{\beq}{\begin{equation}}
\newcommand{\eeq}{\end{equation}}

\numberwithin{equation}{section}
\numberwithin{figure}{section}

\begin{document}
\title{Layer Potential Techniques for the Narrow Escape Problem}

\author{Habib Ammari\thanks{\footnotesize  Department of Mathematics and Applications, Ecole Normale Sup\'erieure,
45 Rue d'Ulm, 75005 Paris, France (habib.ammari@math.cnrs.fr).}
\and Kostis Kalimeris\thanks{\footnotesize Center of Applied
Mathematics, Ecole Polytechnique, 91128 Palaiseau Cedex, France
(kalimeris@cmapx.polytechnique.fr). } \and Hyeonbae
Kang\thanks{\footnotesize Department of Mathematics, Inha
University, Incheon 402-751, Korea (hbkang@inha.ac.kr,
hdlee@inha.ac.kr).} \and Hyundae Lee\footnotemark[3]}

\date{}
\maketitle

\begin{abstract}
The narrow escape problem consists of deriving the asymptotic
expansion of the solution of a drift-diffusion equation with the
Dirichlet boundary condition on a small absorbing part of the
boundary and the Neumann boundary condition on the remaining
reflecting boundaries. Using layer potential techniques, we
rigorously find high-order asymptotic expansions of such
solutions. We explicitly show the nonlinear interaction of many
small absorbing targets. Based on the asymptotic theory for
eigenvalue problems developed in \cite{book}, we also construct
high-order asymptotic formulas for eigenvalues of the Laplace and
the drifted Laplace operators for mixed boundary conditions on
large and small pieces of the boundary.

\end{abstract}

\bigskip

\noindent {\footnotesize Mathematics Subject Classification
(MSC2000): 35B40, 92B05}


\section{Introduction}

An interesting  problem in cellular and molecular biology is to
estimate the mean sojourn time, also called the narrow escape
time,  of a Brownian particle in a bounded domain $\Omega$ before
it escapes through a small absorbing arc $\partial \Omega_a$ on
its boundary $\partial \Omega$. The remaining part of the boundary
$\partial \Omega_r =
\partial \Omega \setminus \overline{\partial \Omega_a}$ is assumed
reflecting for the particle. The small arc often represents a
small target on a cellular membrane. Its physiological role is to
regulate flux, which carries a physiological signal \cite{schuss},
\cite{holcman-rev}.

The narrow escape problem is connected to that of calculating the
solution $u_\ep$ of a mixed Dirichlet-Neumann boundary value
problem in $\Omega$, whose Dirichlet boundary is only $\partial
\Omega_a$ on the otherwise Neumann boundary. The escape time can
be estimated asymptotically in the limit $\ep := |\partial
\Omega_a|/|\partial \Omega| \rightarrow 0$ by computing the
asymptotic expansion of $u_\ep$ as $\ep \rightarrow 0$.

In this paper, we first consider a purely diffuse model. We
provide mathematically rigorous derivations of the first- and
second-order terms in the asymptotic expansion of the solution
$u_\ep$ as $\ep\rightarrow 0$ in the presence of a single or many
small targets. When two or more Dirichlet targets cluster together
they interact nonlinearly. The clustering may affect significantly
the asymptotics. Then we study the problem of eigenvalue changes
due to the small targets.  Finally, accounting for a drift term,
we generalize our results to a mixed Dirichlet-Neumann boundary
value problem for the drift-diffusion equation.

The narrow escape problem of a free Brownian particle (without
drift) through a small target was discussed in \cite{holcman04},
\cite{singer1}, \cite{singer2}, and \cite{singer3}. The method of
\cite{singer1}, \cite{singer2}, \cite{singer3} was generalized in
\cite{singer-schuss} to obtain the leading-order term of the
solution to the corresponding mixed boundary value problem for the
drift-diffusion equation. Matched asymptotics \cite{ward1},
\cite{ward2}, \cite{ward3}, \cite{ward4} yield the expansion of
the principle eigenvalue of the Laplace operator for mixed
boundary conditions on large and small pieces of the boundary. The
effect of clustering of the Dirichlet targets on the first
eigenvalue was analyzed in \cite{holcman08}. The second-order term
in its asymptotic expansion as the size of the target goes to zero
was provided in \cite{singer} by determining the structure of the
boundary singularity of the Neumann function for the Laplacian in
a bounded smooth domain. In all of these papers, the derivations
are quite formal. It is the purpose of this work to provide a
rigorous framework for systematically deriving high-order
asymptotic formulas for the solutions of the diffusion and
drift-diffusion equations with mixed boundary conditions on large
and small pieces of the boundary and the eigenvalues of the
corresponding operators. Our derivations are based on layer
potential techniques and the asymptotic theory for eigenvalue
problems developed in \cite{book}.

While this work was in progress, we found out the existence of the
work \cite{pwpk}. Some of results of this paper, especially those
in Section 3.1, were also obtained in \cite{pwpk}. However, the
method of derivation of asymptotics in \cite{pwpk} is formal and
uses the method of matched asymptotics.

This paper is organized as follows. In Section 2, we formulate the
problem with and without the drift term, and review some facts on
relevant Neumann functions and layer potentials on arcs. Section 3
is to derive (higher-order) asymptotic formula for the solutions
to the narrow escape problem with a single or multiple
(well-separated and closely located) absorbing regions. Section 4
is to deal with the eigenvalue perturbation problem in the
presence of small absorbing region using the same layer potential
techniques. Section 5 is for the narrow escape problem and the
eigenvalue perturbation problem in the
 presence of the force field. The paper ends with a short discussion.

\section{Formulation of the problem}

\subsection{Physical background}
Let $\Omega$ be a bounded simply connected domain in $\R^2$ with $\mathcal{C}^2$-smooth boundary. Suppose that $\partial\Omega$ has two
disjoint parts, the reflecting part $\partial \Omega_r$ and the absorbing part $\partial \Omega_a$, satisfying
$\partial\Omega=\overline{\partial\Omega_a} \cup
\overline{\partial\Omega_r}$. Both $\p\Om_a$ and $\p\Om_r$ consist of finite number of open arcs.
We assume that $\ep:={|\partial\Omega_a|}/2$ is much smaller than $1$ while $|\p\Om|$ is of order $1$. Here and throughout this paper $|\partial \Omega|$ denotes the arc-length of $\partial \Omega$.

Suppose that a Brownian particle is confined to $\Omega$. The probability density function $p_\ep(x,t)$ of
finding the Brownian particle at location $x$ at time $t$ (prior
to its escape) satisfies the Fokker-Planck equation
$$
\frac{\partial p_\ep}{\partial t}(x,t) = \Delta p_\ep(x,t) -
\nabla \cdot (F(x) p_\ep(x,t))
$$
with the initial condition
$$
p_\ep(x,0) = \rho(x),
$$
and the mixed boundary conditions for $t>0$
$$
\left\{ \begin{array}{ll}
p_\ep = 0 \quad & \mbox{on } \partial \Omega_a,\\
\nm \ds \frac{\partial p_\ep}{\partial \nu} - p_\ep F\cdot \nu =0
\quad &\mbox{on } \partial \Omega_r.
\end{array}
\right.
$$
The force field $F$ is given by $F(x)= \nabla \phi(x)$ for a smooth potential
$\phi$. The function $\rho(x)$ is the initial probability
density function; e.g., $\rho(x)= 1/|\Omega|$ for a uniform
distribution or $\rho(x)= \delta_y$, a Dirac mass at $y$,  when
the particle is initially located at $y$.

The function $v_\ep(x) := \int_0^{+\infty} p_\ep(x,t) \,dt$, which
is the mean time the particle spends at $x$ before it escapes
through $\partial \Omega_a$, is the solution to
$$
\left\{\begin{array}{ll} \Delta v_\ep - \nabla \cdot (F v_\ep) = -
\rho \quad &\mbox{in } \Omega,\\
\nm v_\ep = 0 \quad &\mbox{on } {\partial \Omega_a},\\
\nm \ds \frac{\partial v_{\ep}}{\partial \nu}  - v_\ep F\cdot \nu
=0 \quad &\mbox{on }{\partial \Omega_r}.
\end{array}
\right.
$$
The function $w_\ep(x):= v_\ep(x) e^{-\phi(x)}$ is the solution of
the adjoint problem
$$
\left\{\begin{array}{ll} \Delta w_\ep + F \cdot \nabla w_\ep = -
\rho e^{-\phi} \quad &\mbox{in } \Omega,\\
\nm w_\ep = 0 \quad &\mbox{on } {\partial \Omega_a},\\
\nm \ds \frac{\partial w_{\ep}}{\partial \nu} =0 \quad &\mbox{on
}{\partial \Omega_r}.
\end{array}
\right.
$$

Suppose that $\rho(x)= \delta_y$. The function
$$ \ds u_\ep(x)= \frac{\int_\Omega w_\ep(x,y)\, dy}{\int_\Omega e^{-
\phi(y)}\, dy}$$ is the solution to the mixed boundary value
problem
\begin{equation}\label{Poisson}
\begin{cases}
\Delta u_{\ep} + F\cdot \nabla u_\ep =-1 \quad  &\mbox{in }\Omega,\\
u_{\ep}=0 \quad &\mbox{on }{\partial \Omega_a},\\
\nm
\ds \frac{\partial u_{\ep}}{\partial \nu}=0 \quad
&\mbox{on }\partial \Omega_r.
\end{cases}
\end{equation}

The narrow escape problem consists of deriving the asymptotic
expansion of $u_\ep$ as $\ep \rightarrow 0$, from which one can estimate the sojourn time of the Brownian particle. A closely related problem is to construct the asymptotic of the principal eigenvalue
of $-\Delta - F\cdot \nabla$ in $\Omega$ with the mixed
Dirichlet-Neumann boundary conditions for small $\ep$.

The goal of this paper is threefold: (i) To find the asymptotic
behavior of the solution $\ue$ to \eqnref{Poisson} with $F=0$
(free Brownian motion) in the presence of one or multiple small
targets as $\ep \rightarrow 0$; (ii) To find the asymptotics of
the eigenvalues of $-\Delta$ in $\Omega$ with the Dirichlet
boundary condition on $\partial \Omega_a$ and the Neumann boundary
condition on $\partial \Omega_r$; (iii) To generalize these
asymptotics to the case with nonzero drift ($F\neq 0$).

We note that if the force field $F=0$, then  by integrating
\eqnref{Poisson} over $\Omega$ the following compatibility
condition holds:
\begin{equation}\label{C_1}
\int_{\partial \Omega_a} \frac{\partial \ue}{\partial
\nu}\, d\sigma=-|\Omega|.
\end{equation}

\subsection{Neumann functions}

Let $N(x,z)$ be the Neumann function for $-\Delta$ in $\Omega$
corresponding to a Dirac mass at $z \in \Om$. That is, $N$ is the solution
to
\begin{equation}\label{Neumann}
\begin{cases}
&\Delta_z N(x,z)=-\delta_x, \ \ x,z\in \Omega,\\
& \displaystyle\frac{\partial N}{\partial \nu_z}\bigg|_{z \in \partial \Omega}
=-\frac{1}{|\partial\Omega|}, \ \ \ \int_{\partial \Omega}N(x,z) d\sigma(z)=0, \ \ x\in\Omega.
\end{cases}
\end{equation}
The Neumann function is symmetric in its arguments;
$N(x,z)=N(z,x)$ for $x\neq z \in \Omega$. It furthermore takes the
form
\begin{equation}
\label{Neumann smooth}N(x,z)=-\frac{1}{2\pi} \ln|x-z| +
R_\Omega(x,z), \quad x, z\in\Omega,
\end{equation}
where $R_{\Omega}(\cdot,z)$ belongs to $H^{{3}/{2}}(\Omega)$, the
standard Sobolev space of order $3/2$, for any $z\in\Omega$ and
solves
\begin{equation}
\left\{
\begin{aligned}
\label{R_w smooth}
&\Delta_x R_\Omega(x,z)=0, \ \ x\in \Omega,\\
& \frac{\partial R_\Omega}{\partial \nu_x}\bigg|_{x \in \partial
\Omega}=-\frac{1}{|\partial\Omega|}+\frac{1}{2\pi}\frac{\langle
x-z,\nu_x \rangle}{|x-z|^2}, \ \ x\in\partial\Omega.
\end{aligned}
\right.
\end{equation}

Suppose $z\in\partial\Omega$. Since $\partial \Omega$ is smooth,
only a half of any sufficiently small disk about $z$
is contained in $\Omega$, and hence the singularity of
$N(x,z)$ is $-(1/\pi) \ln |x-z|$. Therefore, $N(x,z)$ for $z \in \partial \Omega$, which we
denote by $N_{\partial\Omega}$, can be
written as
\begin{equation}
\label{Neumann smooth
boundary}N_{\partial\Omega}(x,z)=-\frac{1}{\pi} \ln|x-z| +
R_{\partial\Omega}(x,z), \ \ x\in \Omega,\ z\in\partial\Omega,
\end{equation}
where $R_{\partial\Omega}(\cdot,z)$ solves the problem
\begin{equation}\label{R_tw smooth}
\left\{
\begin{aligned}
&\Delta_x R_{\partial\Omega}(x,z)=0, \ \ x\in \Omega,\\
& \frac{\partial R_{\partial\Omega}}{\partial \nu_x}\bigg|_{x \in
\partial
\Omega}=-\frac{1}{|\partial\Omega|}+\frac{1}{\pi}\frac{\langle
x-z,\nu_x\rangle}{|x-z|^2}, \ \ x\in\partial\Omega.
\end{aligned}
\right.
\end{equation}

Note that the Neumann data in \eqnref{R_tw smooth} is bounded
on $\p\Om$ uniformly in $z \in \p\Om$ since $\p\Om$ is $\mathcal{C}^2$-smooth,
and hence $R_{\partial\Omega}(\cdot,z)$ belongs to
$H^{{3}/{2}}(\Omega)$ uniformly in $z\in\partial\Omega$.

Let
 \beq\label{definegx}
 g(x) := \int_\Omega  N(x,z)dz, \quad x \in \Om.
 \eeq
Then $\Delta g =-1$ in $\Om$ and $\pd{g}{\nu}= -|\Om|/|\p\Om|$ on $\p\Om$.
Therefore one may use Green's formula and the mixed boundary conditions in
\eqnref{Poisson} with $F=0$ to have
\begin{equation}
\label{Int representation}\ue(x)=g(x) +
\int_{\partial \Omega_a} N_{\p\Om}(x,z) \frac{\partial \ue(z)}{\partial
\nu_z} d\sigma(z)+\Ce, \quad x \in \Omega,
\end{equation}
where
 \beq
 \Ce=\frac{1}{|\partial \Omega|} \int_{\partial
 \Omega_r}\ue(z)d\sigma(z).
 \eeq
In view of \eqnref{Neumann
smooth boundary}, \eqnref{Int representation} becomes
\begin{equation}\label{Int representation smooth}
\ue(x)=g(x)-\frac{1}{\pi} \int_{\partial
\Omega_a}\ln|x-z|\frac{\partial \ue(z)}{\partial \nu} d\sigma(z)+
\int_{\partial \Omega_a}R_{\partial\Omega}(x,z)\frac{\partial
\ue(z)}{\partial \nu} d\sigma(z)+\Ce.
\end{equation}

If $\Om$ is the unit disk, then
 $$
 \frac{\langle x-z,\nu_x \rangle}{|x-z|^2}=\frac{\langle x-z,x
 \rangle}{|x|^2+|z|^2-2x\cdot z}=\frac{1-x\cdot z}{2-2x\cdot z}=\frac{1}{2}
 $$
for $x\in\partial\Omega$. Therefore,
$\frac{\partial R_{\p\Omega}}{\partial \nu_x} |_{x \in \partial
\Omega}=0$, and hence $R_{\partial\Omega}(x,z)=\mbox{const}$. Because of the condition $\int_{\p\Om} N(x,z) d\sigma(z)=0$, we have $R_{\partial\Omega}(x,z)=0$ for all $x\in\Omega$ and $z\in \partial\Omega$, and hence
 \beq\label{diskGreen}
 N_{\partial\Omega}(x,z)=-\frac{1}{\pi} \ln|x-z|, \ \ x\in \Omega,\ z\in\partial\Omega.
 \eeq
We also have
 \beq\label{cepdisk}
 g(x)= \int_\Omega N(x,z)dz=\frac{1}{4}(1-|x|^2),
 \eeq
which can be seen using Green's formula. In fact, we have
\begin{equation*}
\int_\Omega N(x,z) \Delta |z|^2 - \Delta_z N(x,z) |z|^2 dz =
\int_{\partial \Omega}  \bigg(2N(x,z)|z| - \frac{\partial
N(x,z)}{\partial \nu_z} |z|^2 \bigg)\, d\sigma(z).
\end{equation*}

\subsection{Single-layer potential on an arc}

Our aim is to derive an asymptotic expansion  of ${\partial
\ue(z)}/{\partial \nu}$ on $\partial \Omega_a$ for $\ep$ small
enough. We write this density function as the solution of an integral
equation \eqnref{Int representation smooth} on $\partial \Omega_a$. For doing so, we introduce
two Hilbert spaces: For $\ep >0$, define
 \beq
 \bX_\ep : =  \big\{ \varphi: \int_{-\ep}^\ep \sqrt{\ep^2 - x^2}\,
 |\varphi(x)|^2 \, dx < +\infty \big\}
 \eeq
with the norm
 \[
 \| \varphi \|_{\mathbf{X}_\ep} = \big( \int_{-\ep}^\ep \sqrt{\ep^2 - x^2}\,
 |\varphi(x)|^2 \, dx  \big)^{1/2}.
 \]
Then one immediately has
 \beq
 \| \varphi \|_{\mathbf{X}_\ep} = \| \widetilde{\varphi} \|_{\mathbf{X}_1},
 \eeq
where $\tilde{\varphi}(x):= \ep \varphi(\ep x)$. We also define
 \beq
 \mathbf{Y}_\ep := \big\{  \psi \in \mathcal{C}^0 \left(\,\left[\, - \ep,
 \ep \,\right]\,\right): \psi^{\prime} \in \mathbf{X}_\ep \big\},
 \eeq
with the norm
\[\| \psi \|_{\mathbf{Y}_\ep} = \big( \|\psi\|^2_{\mathbf{X}_\ep} + \|
\psi^\prime\|^2_{\mathbf{X}_\ep} \big)^{1/2},\]
 where $\psi^{\prime}$ is the derivative of $\psi$ in the sense of distribution.


We also recall the following lemma; see \cite[Chapter 5]{book} and
the references therein.
\begin{lemma}\label{invers1}
The integral operator ${\mathcal L}: {\mathbf X}_1 \mapsto {\mathbf
Y}_1$ defined by
 \beq\label{operatorL}
 \mathcal{L}[\varphi](x)=\int_{-1}^1 \ln |x-y|\varphi(y)dy
 \eeq
is invertible. For a given function $\psi \in {\mathbf Y}_1$,
$\mathcal{L}^{-1}[\psi] \in {\mathbf X}_1$ is given by
 \begin{equation} \label{invcarleman}
 \mathcal{L}^{-1}[\psi](x) \,=\,  - \frac{1}{\pi^2
 \sqrt{1 - x^2}}\, \int_{-1}^{1}
 \frac{\sqrt{1 - y^2} \psi^{\prime}(y)}{x - y}\; dy -
 \frac{a(\psi)}{ \pi (\ln 2)\, \sqrt{1
 - x^2} }
 \end{equation}
for $x \in (-1, 1 )$,
where the constant $a(\psi)$ is defined by
 \begin{equation} \label{defaf}
 a(\psi)=\psi(x) +
 {\mathcal L} \bigg[\frac{1}{\pi^2 \sqrt{1 - y^2}}\,
 \int_{- 1}^{1} \frac{\sqrt{1 - z^2}
 \psi^{\prime}(z)}{y - z}\; dz \bigg] (x).
 \end{equation}
\end{lemma}

Note that, in view of \eqnref{invcarleman} and \eqnref{defaf}, we
have
 \begin{equation} \label{invcarleman1}
 \mathcal{L}^{-1}[1](x) = - \frac{1}{ \pi (\ln {2})  \sqrt{1 - x^2}}.
 \end{equation}

\section{Narrow escape of a free Brownian particle}

We now consider the narrow escape problem in the presence of a single or multiple (clustered or well-separated) absorbing arcs.

\subsection{A single small target}

Let  $x(t): [-\ep,\ep]\to \R^2$ be the arclength parametrization
of $\partial \Omega_a$; namely, $x$ is a $\mathcal{C}^2$-function
satisfying $|x'(t)|= 1$ for all $t \in [-\ep, \ep]$ and
 \begin{equation} \label{paramet} \partial
 \Omega_a=\{x(t):t\in[-\ep,\ep]\}.\end{equation}

Since $u_\ep=0$ on $\p\Om_a$, it follows from \eqnref{Int representation smooth} that
\begin{equation}\label{int_eq203}
f(t)-\frac{1}{\pi} \int_{-\ep}^\ep\ln|x(t)-x(s)|\fe(s) ds + \int_{-\ep}^\ep r(t,s) \fe(s) ds + \Ce=0,
\end{equation}
where
 $$
 f(t):=g(x(t)),\quad
 \fe(t):=\frac{\partial \ue(x(t))}{\partial \nu},
 \quad  r(t,s):=R_{\partial\Omega}(x(t),x(s)).
 $$
By a simple change of variables, we obtain that
\begin{equation} \label{eq3.3}
\frac{1}{\pi} \int_{-1}^1\ln|x(\ep t)-x(\ep s)|\fte(s) ds -
\int_{-1}^1 r(\ep t, \ep s) \fte(s)ds = f(\ep t) + \Ce,
\end{equation}
with $\fte(t)=\ep\fe(\ep t).$

Since the compatibility condition \eqnref{C_1} reads
 \begin{equation}\label{C_11}
 \int_{-1}^1 \fte(s)ds=-|\Omega|,
 \end{equation}
we may rewrite \eqnref{eq3.3} as
\begin{equation}\label{int_eq103}
\frac{1}{\pi}(\mathcal{L}+\ep\mathcal{L}_1)[\fte](t) -\frac{|\Omega|\ln \ep}{\pi}+  r(0,0)|\Omega| = f(\ep t)+ \Ce,
\end{equation}
where the operator $\mathcal{L}_1$ is given by
 $$
 \mathcal{L}_1[\varphi]:=\frac{1}{\ep}\int^1_{-1} \left(\ln\frac{|x(\ep t)-x(\ep s)|}{\ep|t-s|} + \pi r(0,0)- \pi r(\ep t,\ep s)\right) \varphi(s) ds.
 $$
Moreover, since
$$ |x(\ep t)-x(\ep s)|=\ep|t-s|(1+O(\ep)),$$
one can see that $\mathcal{L}_1:\mathbf{X}_1 \mapsto
\mathbf{Y}_1$ is bounded independently of $\ep$.
Note that $f$ is $\mathcal{C}^1$ on $\overline{\Om}$ and hence we may rewrite \eqnref{int_eq103} as
\begin{equation}
(\mathcal{L}+\ep\mathcal{L}_1)[\fte](t)= {|\Omega|\ln \ep}- \pi r(0,0)|\Omega| + \pi\Ce+\pi f(0)+ O(\ep).
\end{equation}
It then follows from \eqnref{invcarleman1} that
\begin{equation}
\label{fte
smooth}\fte(t)=\left[\frac{|\Omega|\ln\ep}{\pi}-r(0,0)|\Omega| +
\Ce+ f(0)\right] \frac{1}{(\ln\frac{1}{2})\, \sqrt{1-t^2}}+
O(\ep),
\end{equation}
where the remainder $O(\ep)$ is with respect to the norm
$\|\cdot\|_{\mathbf{X}_1}$.

Now, plugging \eqnref{fte smooth} into \eqnref{C_11}  we arrive at
 \begin{equation}\label{C_111}\Ce=\frac{|\Omega|}{\pi}\ln\frac{2}{\ep}+r(0,0)|\Omega|- f(0)+O(\ep)\end{equation}
 and, consequently,
\begin{equation}
\label{fte smooth1}\fte(t)= -\frac{|\Omega|}{\pi}\frac{1}{\sqrt{1-t^2}} + O(\ep),
\end{equation}
where $O(\ep)$ is with respect to $ \|\cdot \|_{\mathbf{X}_1}$. Therefore, we obtain
\begin{equation}\label{fesmooth}
\fe(t)= -\frac{|\Omega|}{\pi}\frac{1}{\sqrt{\ep^2-t^2}} + O(\ep),
\end{equation}
where $O(\ep)$ is  with respect to $ \|\cdot \|_{\mathbf{X}_\ep}$.

Finally, combining \eqnref{Int representation smooth} and
\eqnref{fesmooth} yields
\begin{align*}
\label{solution smooth}
\ue(x)&=g(x)-\frac{|\Omega|}{\pi} \int_{-\ep}^{\ep} \frac{N_{\partial\Omega}(x,x(t))}{\sqrt{\ep^2-t^2}} dt +\Ce+O(\ep)\\
&=g(x) -|\Omega| N_{\partial\Omega}(x,x^*)+\Ce+O(\ep)
\end{align*}
for $x\in \Omega$ provided that $\mbox{\rm
dist}(x,\partial\Omega_a) \ge c$ for some constant $c>0$. Thus we have the following theorem.

\begin{theorem} \label{themF0}
Suppose that $\p\Omega_a$ is an arc of center $x^*$ and length $2\ep$. Then
the following asymptotic expansion of $u_{\ep}$ holds
\begin{equation}
\label{solution smooth 1}
\ue(x)=\frac{|\Omega|}{\pi} \ln\frac{2}{\ep}+ \Phi_\Omega(x,x^*) +O(\ep),
\end{equation}
where
 \beq\label{surfaceGreen}
 \Phi_\Omega(x,x^*)=\int_\Omega N(x,z)dz - |\Omega| N_{\partial\Omega}(x,x^*)-\int_\Omega N(x^*,z)dz+|\Omega|R_{\partial\Omega}(x^*,x^*).
 \eeq
The remainder $O(\ep)$ is uniform in $x\in \Omega$ satisfying
$\mbox{\rm dist}(x,\partial\Omega_a) \ge c$ for some constant $c>0$. Moreover, if $x(t)$,
$-\ep<t<\ep$, is the arclength parametrization of $\partial
\Omega_a$, then
$$ \frac{\partial u_{\ep}}{\partial \nu}(x(t)) =-\frac{|\Omega|}{\pi}\frac{1}{\sqrt{\ep^2-t^2}} + O(\ep),$$
where $O(\ep)$ is with respect to $\|\cdot\|_{\bX_\ep}$.
\end{theorem}

The leading order term in \eqnref{solution smooth 1} was obtained
in \cite{singer1} and the formula \eqnref{solution smooth 1} was
obtained in \cite{pwpk} in a formal way (see (2.14) in that
paper). We emphasize that the derivation of this paper is
mathematically rigorous. Moreover, the method of derivation is
quite different from that of \cite{pwpk}; It is based on the layer
potential technique. We note that the corrector function
$\Phi_\Omega(x,x^*)$ solves the following problem
\begin{equation} \label{phiom}
\left\{\begin{array}{ll} \Delta_x \Phi_\Omega(x,x^*) = -1 \quad
& \mbox{in } \Omega,\\
\nm
\ds \frac{\partial \Phi_\Omega}{\partial \nu_x} = - |\Omega|
\delta_{x^*} \quad & \mbox{on } \partial \Omega,
\end{array}
\right.
\end{equation}
which shows that \eqnref{solution smooth 1} is nothing else than a
dipole-type approximation.

If $\Omega$ is the unit disk centered at $0$, one can easily see from \eqnref{diskGreen} and \eqnref{cepdisk} that
$$\Phi_\Omega(x,x^*)=\ln|x-x^*| + \frac{1}{4}(1-|x|^2).$$
Thus we have the following corollary.
\begin{cor} Suppose that $\Omega$ is the unit disk and $\partial\Omega_a$ is an arc of center $x^*$ and length $2\ep$.
Then the solution $u_{\ep}$ is given asymptotically as
\begin{equation}\label{uedisk}
\ue(x)=\ln\frac{2}{\ep}+\frac{1}{4}(1-|x|^2)+\ln|x-x^*| +O(\ep),
\end{equation}
where the remainder $O(\ep)$ is uniform in $x\in \Omega$ satisfying
$\mbox{\rm dist}(x,\partial\Omega_a) \ge c$ for some constant $c$.
\end{cor}

The formula \eqnref{uedisk} was also obtained in \cite{singer2} (and \cite{pwpk}).


\subsection{Two well-separated small targets}

Let $\Omega$ be a unit disk centered at $0$. Suppose
that $\p\Om_a$ consists of two parts, say $\p\Om_1$ and $\p\Om_2$, satisfying
$\overline{\partial\Omega_1} \cap \overline{\partial\Omega_2}=\emptyset$.

Let $u$ be the solution to \eqnref{Poisson} with $F=0$. In view of \eqnref{Int representation smooth},
$u$ is now written as
\begin{equation} \label{Int representation_2}
u(x)= g(x) + \int_{\partial \Omega_1} N_{\p\Om}(x,z)\frac{\partial u(z)}{\partial \nu} d\sigma(z) +  \int_{\partial \Omega_2}N_{\p\Om}(x,z) \frac{\partial u(z)}{\partial \nu} d\sigma(z) +\Ce.
\end{equation}

Suppose that $\partial \Omega_i$, $i=1,2$, are parameterized by
 $$x(t)=(\cos t,\sin t),\ t\in(s_i-\ep_i,s_i+\ep_i).$$
Then the (Dirichlet) boundary conditions on $\p\Om_a$ yield the following integral equations
\begin{equation}\label{integral eq1_2targets}
0= \int_{s_1-\ep_1}^{s_1+\ep_1}N_{\p\Om} (x(s),x(t))\phi(t) dt + \int_{s_2-\ep_2}^{s_2+\ep_2}N_{\p\Om} (x(s),x(t))\phi(t) dt + \Ce
\end{equation}
for $s\in(s_i-\ep_i,s_i+\ep_i)$ where $\phi(t):=\pd{u}{\nu}(x(t))$.

Recall that  $N_{\p\Om} (x,z)=-\frac{1}{\pi}\ln|x-z|,\ x\in
\partial\Omega,\ z\in \partial\Omega$ for  the unit disk. If we
make changes of variables $s\mapsto s_i+\ep_i s$ and $t\mapsto s_i+\ep_i t$
for $s,t\in(s_i-\ep_i,s_i+\ep_i)$, the integral equations \eqnref{integral eq1_2targets} become
\begin{equation}\label{integral eq3_2targets}
\begin{cases}
&\ds\int_{-1}^{1} \ln{|x(s_1+\ep_1 s)-x(s_1+\ep_1 t)|} \, {\fo (t)} \,dt \\
&\qquad \qquad \ds +   \int_{-1}^{1} \ln{|x(s_1+\ep_1 s)-x(s_2+\ep_2 t)|} \, {\ft (t)} \,dt =\pi \Ce,\\
&\ds\int_{-1}^{1} \ln{|x(s_2+\ep_2 s)-x(s_1+\ep_1 t)|} \, {\fo (t)} \, dt \\
& \qquad \qquad \ds +   \int_{-1}^{1} \ln{|x(s_2+\ep_2
s)-x(s_2+\ep_2 t)|} \, {\ft (t)} \, dt =\pi \Ce,
\end{cases}
\end{equation}
where  $\varphi_i(t)=\ep_i\phi(s_i+\ep_i t)$ for $ t\in(-1,1)$.

Let
 \beq
 C_i:=\frac{1}{\pi}\int^1_{-1} \varphi_i (t)dt, \quad i=1,2.
 \eeq
Because of \eqnref{C_1}, we have
\begin{equation}\label{C_j b2 targets}
C_1+C_2=-1.
\end{equation}
One can easily check that
\begin{equation}
\begin{cases}
\ln{|x(s_i+\ep_i s)-x(s_i+\ep_i t)|}=\ln\ep_i|s-t|+ O(\ep_i),\\
\nm
\ln{|x(s_i+\ep_i s)-x(s_j+\ep_j t)|}=\ln|x(s_1)-x(s_2)|+O\left(\frac{\sqrt{\ep_1^2+\ep_2^2}}{|x(s_1)-x(s_2)|}\right),~i\ne j.
\end{cases}
\end{equation}
Therefore, the system \eqnref{integral eq3_2targets} can be written in
the form
\begin{equation}
\label{integral eq4_2targets_2}
(\mathcal{A}+\mathcal{B})\begin{pmatrix} \fo \\ \ft \end{pmatrix}=
\pi\begin{pmatrix}  \Ce-C_1\ln\ep_1-C_2\ln|x(s_1)-x(s_2)|\\  \Ce-C_1\ln|x(s_1)-x(s_2)|-C_2\ln\ep_2\end{pmatrix},
\end{equation}
where the operator $\mathcal{A}:\mathbf{X}_1 \times \mathbf{X}_1 \mapsto \mathbf{Y}_1 \times \mathbf{Y}_1$ is given by
$$ \mathcal{A}:=\begin{pmatrix} \mathcal{L} & 0 \\ 0 & \mathcal{L}\end{pmatrix}, $$
and
$\mathcal{B}:\mathbf{X}\times\mathbf{X}\mapsto\mathbf{Y}\times\mathbf{Y}$ satisfies
 $$
 \mathcal{B}=\begin{pmatrix} O(\ep_1) & O\left(\frac{\sqrt{\ep_1^2+\ep_2^2}}{|x(s_1)-x(s_2)|}\right) \\
 O\left(\frac{\sqrt{\ep_1^2+\ep_2^2}}{|x(s_1)-x(s_2)|}\right) & O(\ep_2)\end{pmatrix}
 $$
in the operator norm. Here $\mathcal{L}$ is given by \eqnref{operatorL}.

Let us now assume that $|x(s_1)-x(s_2)|\ge c$ for some $c>0$, in other words, the two
targets are well-separated. It then follows from \eqnref{invcarleman1} that
\begin{equation}\label{integral eq6_2targets}
\begin{split}
&\fo(s) = \frac{\Ce-C_1\ln\ep_1-C_2\ln|x(s_1)-x(s_2)|}{\ln\frac{1}{2}} \left(\frac{1}{\sqrt{1-s^2}}
+ O(\sqrt{\ep_1^2+\ep_2^2})\right),\\
&\ft(s) = \frac{\Ce-C_1\ln|x(s_1)-x(s_2)|-C_2\ln\ep_2}{\ln\frac{1}{2}} \left( \frac{1}{\sqrt{1-s^2}} + O(\sqrt{\ep_1^2+\ep_2^2})\right),
\end{split}
\end{equation}
where $O(\sqrt{\ep_1^2+\ep_2^2})$ is with respect to
$\|\cdot\|_{\mathbf{X}_1}$.

Integrating \eqnref{integral eq6_2targets} over $(-1,1)$, we obtain a system of two equations for $C_1,C_2$ and $\Ce$
\begin{equation}
\label{C_j a2 targets}
\begin{cases}
\ds C_1\ln\frac{\ep_1}{2}+ C_2\ln |x(s_1)-x(s_2)| - \Ce = O(\sqrt{\ep_1^2+\ep_2^2}),\\
\nm
\ds C_1\ln |x(s_1)-x(s_2)| + C_2 \ln\frac{\ep_2}{2} -\Ce = O(\sqrt{\ep_1^2+\ep_2^2}),
\end{cases}
\end{equation}
which together with \eqnref{C_j b2 targets} yields
 \begin{align*}
 C_1 &=-\frac{\ln\frac{\ep_2}{2|x(s_1)-x(s_2)|}}{\ln\frac{\ep_1\ep_2}{4|x(s_1)-x(s_2)|^2}}+O(\sqrt{\ep_1^2+\ep_2^2}),\\ C_2 &=-\frac{\ln\frac{\ep_1}{2|x(s_1)-x(s_2)|}}{\ln\frac{\ep_1\ep_2}{4|x(s_1)-x(s_2)|^2}}+O(\sqrt{\ep_1^2+\ep_2^2}), \\
 \Ce &=\frac{(\ln |x(s_1)-x(s_2)|)^2 - \ln \frac{\ep_1}{2}\ln \frac{\ep_2}{2}}{\ln\frac{\ep_1\ep_2}
 {4|x(s_1)-x(s_2)|^2}}+O(\sqrt{\ep_1^2+\ep_2^2}).
 \end{align*}
Thus we have
\begin{equation*}
\begin{split}
&\fo(s) =-\frac{\ln\frac{\ep_2}{2|x(s_1)-x(s_2)|}}
{\ln\frac{\ep_1\ep_2}{4|x(s_1)-x(s_2)|^2}} \frac{1}{\sqrt{1-s^2}}  + O(\sqrt{\ep_1^2+\ep_2^2}),\\
&\ft(s) =-\frac{\ln\frac{\ep_1}{2|x(s_1)-x(s_2)|}}{\ln\frac{\ep_1\ep_2}{4|x(s_1)-x(s_2)|^2}} \frac{1}{\sqrt{1-s^2}} + O(\sqrt{\ep_1^2+\ep_2^2}),
\end{split}
\end{equation*}
where $O(\sqrt{\ep_1^2+\ep_2^2})$ is with respect to
$\|\cdot\|_{\mathbf{X}_1}$. By scaling back, the
following expansion of $\phi(s)$ holds as $\ep \to 0$:
\begin{equation*}
\label{integral eq7_2targets}
\phi(s)=\begin{cases}
 -\ds \frac{\ln\frac{\ep_2}{2|x(s_1)-x(s_2)|}}{\ln\frac{\ep_1\ep_2}{4|x(s_1)-x(s_2)|^2}} & \ds \frac{1}{\sqrt{\ep_1^2-(s-s_1)^2}} + O(\sqrt{\ep_1^2+\ep_2^2}), \quad s\in(s_1-\ep_1,s_1+\ep_1),\\
 \nm
 \nm
 -\ds \frac{\ln\frac{\ep_1}{2|x(s_1)-x(s_2)|}}{\ln\frac{\ep_1\ep_2}{4|x(s_1)-x(s_2)|^2}}   & \ds\frac{1}{\sqrt{\ep_2^2-(s-s_2)^2}} + O(\sqrt{\ep_1^2+\ep_2^2}),\quad s\in(s_2-\ep_2,s_2+\ep_2).
\end{cases}
\end{equation*}
By substituting this formula to \eqnref{Int representation_2}, we obtain the following theorem.
\begin{theorem}
Let  $\partial \Omega_i$, $i=1,2$ be parameterized by
 $$x(t)=(\cos t,\sin t),\ t\in(s_i-\ep_i,s_i+\ep_i).$$
Assume that there is a constant $c_0>0$ such that $|x(s_1)-x(s_2)| \ge c_0$.
Then the solution $u$ to \eqnref{Poisson} is given asymptotically
by
\begin{equation} \label{two-target-sep}
\begin{split}
u(x)&=\frac{(\ln |x(s_1)-x(s_2)|)^2 -
\ln \frac{\ep_1}{2}\ln \frac{\ep_2}{2}}{\ln\frac{\ep_1\ep_2}{4|x(s_1)-x(s_2)|^2}}+
\frac{1}{4}(1-|x|^2)\\
&+\frac{\ln\frac{\ep_2}{2|x(s_1)-x(s_2)|}}{\ln\frac{\ep_1\ep_2}{4|x(s_1)-x(s_2)|^2}} \ln|x-x(s_1)|+\frac{\ln\frac{\ep_1}{2|x(s_1)-x(s_2)|}}{\ln\frac{\ep_1\ep_2}{4|x(s_1)-x(s_2)|^2}}\ln|x-x(s_2)|+ O(\sqrt{\ep_1^2+\ep_2^2}),
\end{split}
\end{equation}
where $O(\sqrt{\ep_1^2+\ep_2^2})$ is uniform in $x\in \Omega$ satisfying
$\mbox{\rm dist}(x,\partial\Omega_1\cup\partial\Omega_2) \ge c$ for some constant $c>0$.
\end{theorem}

If $\ep:=\ep_1=\ep_2$, then \eqnref{two-target-sep} takes the following simpler form:
\begin{equation} \label{fsum}
\begin{split}
u(x)=\frac{1}{2}\ln \frac{2}{\ep}
+\frac{1}{4}(1-|x|^2) +\frac{1}{2} \ln \frac{|x-x(s_1)||x-x(s_2)|}{|x(s_1)-x(s_2)|} + O(\ep).
\end{split}
\end{equation}
Formula \eqnref{fsum} shows that the leading order term of the escape time when there are two targets of the same size is reduced to half of that when there is a single target, which is quite natural. The second-order term (dipole type) shows that the escape time is reduced to half of the sum of two dipoles increased by the relative position of the targets.

\subsection{Two clustered small targets on the unit circle}
We now consider the case when there are two clustered targets. We suppose that
\begin{itemize}
\item $\ep_1=\ep_2(=:\ep),$
\item $s_2-s_1=d\ep$, $d>2$.
\end{itemize}
We still assume that $\Omega$ is a unit disk centered at $0$.

Let
\begin{align*}
&\mathcal{L}_1[\phi](s):=\int^1_{-1} \ln|d+s-t| \, \phi(t) \, dt,\\
&\mathcal{L}_2[\phi](s):=\int^1_{-1} \ln|d-s+t| \, \phi(t) \, dt.
\end{align*}
The system of equations \eqnref{integral eq3_2targets} now takes
the form
\begin{equation}
\label{integral eq4_2targets_6}
(\mathcal{A}+\mathcal{B})\begin{pmatrix} \fo \\ \ft \end{pmatrix}=
\pi\begin{pmatrix}  \Ce-C_1\ln\ep-C_2\ln\ep\\  \Ce-C_1\ln\ep-C_2\ln\ep\end{pmatrix},
\end{equation}
where  $\mathcal{A}, ~\mathcal{B}:\mathbf{X}\times\mathbf{X}\mapsto\mathbf{Y}\times\mathbf{Y}$ are given by
$$ \mathcal{A}:=\begin{pmatrix} \mathcal{L} & \mathcal{L}_2 \\ \mathcal{L}_1 & \mathcal{L}\end{pmatrix}$$
and $\mathcal{B}= O(\ep)$ in the operator norm.

The following invertibility result follows from \cite{saranen}.
\begin{lemma} \label{leminv}
The operator
$\mathcal{A}:\mathbf{X}\times\mathbf{X}\mapsto\mathbf{Y}\times\mathbf{Y}$
is invertible. \end{lemma}

It follows from Lemma \ref{leminv} and \eqnref{C_j b2 targets}
that
\begin{equation}\label{int_eq_1938}
\begin{pmatrix} \fo \\ \ft \end{pmatrix}=\pi(\Ce+\ln\ep)\left(\mathcal{A}^{-1}\begin{pmatrix} 1 \\ 1 \end{pmatrix}+O(\ep)\right).
\end{equation}

Define $\alpha_j=\alpha_j(d)$, $j=1,2$, by
 \beq
 \begin{pmatrix} \alpha_1 \\ \alpha_2
 \end{pmatrix} := \int^1_{-1} \mathcal{A}^{-1} \begin{pmatrix} 1 \\ 1 \end{pmatrix} (t) dt.
 \eeq
Then, we have $\alpha_1=\alpha_2(=\alpha)$. In fact, if $(\vp_1,
\vp_2)^T$, where $T$ denotes the transpose, is the solution of
 \beq\label{Aoneone}
 \mathcal{A} \begin{pmatrix} \vp_1 \\ \vp_2 \end{pmatrix} = \begin{pmatrix} 1 \\ 1 \end{pmatrix},
 \eeq
so is $(\vp_2(-t), \vp_1(-t))^T$, and hence we have
 \beq\label{vpone}
 \vp_1(t)=\vp_2(-t).
 \eeq

Integrating \eqnref{int_eq_1938} over $\partial \Omega_a$ yields
$$ \Ce= -\ln\ep - \frac{1}{2\alpha(d)} +O(\ep),$$ which combined
with \eqnref{Int representation_2} and \eqnref{int_eq_1938},
gives the following result.
\begin{theorem}\label{adjacent_target}
Let  $\partial \Omega_i$, $i=1,2$ be parameterized by
 $$x(t)=(\cos t,\sin t),\ t\in(s_i-\ep,s_i+\ep)$$
 with $|s_1-s_2|=d\ep$. Then the solution $u$ to \eqnref{Poisson} is given asymptotically by
\begin{equation} \label{twocluster}
\begin{split}
u(x)=\ln \frac{1}{\ep} - \frac{1}{2\alpha(d)} +\frac{1}{4}(1-|x|^2)
+\frac{1}{2} \ln|x-x(s_1)|+\frac{1}{2} \ln|x-x(s_2)|+ O(\ep),
\end{split}
\end{equation}
where $O(\ep)$ is uniform in $x\in \Omega$ satisfying $\mbox{\rm
dist}(x,\partial\Omega_1\cup\partial\Omega_2) \ge c$ for some $c>0$.
\end{theorem}

In view of \eqnref{vpone}, it follows from \eqnref{Aoneone} that
 $$
 \mathcal{L}[\vp_1](s) + \mathcal{L}_2[\vp_1(-\cdot)](s)=1,
 $$
or
 $$
 \int^1_{-1} [\ln|s-t|- \ln|d-s-t|] \, \vp_1(t) \, dt=1.
 $$
By taking $s=1$ one can see that $\alpha(d) \to \infty$ as $d \to
2+$, and hence the solution converges to the one with an absorbing
arc of length $4\ep$ (modulo $O(\ep)$).

It is interesting to compare the formula \eqnref{twocluster} with
\eqnref{fsum}. The nonlinear interaction between the two small
targets is described by the term $\ln \frac{1}{\ep}-
\frac{1}{2\alpha(d)}$ in \eqnref{twocluster} while it is described
by the term $-\frac{1}{2} \ln {|x(s_1)-x(s_2)|}$ in \eqnref{fsum}.



\subsection{Multiple small targets on the unit circle}
We now extend the above analysis to larger clusters. We consider
clusters $\partial \Omega_i$, $i=1,\ldots,n$, parameterized by
 $$x(t)=(\cos t,\sin t),\ t\in(s_i-\ep,s_i+\ep).$$
We assume that each $ \partial \Omega_i$ is separated from its
neighbors by a distance comparable to $\ep$.

Let $d_{ij}:=\frac{s_i-s_j}{\ep}$. Define the operators
$\mathcal{L}_{ij}$ and $\mathcal{A}$ by
$$ \mathcal{L}_{ij}[\phi](s):=\int^1_{-1} \ln|d_{ij}+s-t|\phi(t)dt$$
and
$$\mathcal{A}:=(\mathcal{L}_{ij})_{n\times n}.$$

Let
$$\begin{pmatrix} \alpha_1\\ \vdots\\ \alpha_n\end{pmatrix}
=\int^1_{-1} \mathcal{A}^{-1} \begin{pmatrix} 1\\ \vdots\\ 1\end{pmatrix}dt.$$

Following the same lines as in the proof of Theorem
\ref{adjacent_target}, we can easily prove the following result.

\begin{theorem}
Let  $\partial \Omega_i$, $i=1,\ldots,n$ be parameterized by
 $$x(t)=(\cos t,\sin t),\ t\in(s_i-\ep,s_i+\ep)$$
 with $s_i-s_j=d_{ij}\ep$. Then the solution to \eqnref{Poisson} is given asymptotically by
\begin{equation}
\begin{split}
u(x)=-\ln\ep - \frac{1}{\alpha_1+\cdots+\alpha_n} +\frac{1}{n}\sum_{i=1}^n \Phi_\Omega(x,x(s_i))+ O(\ep),
\end{split}
\end{equation}
where the remainder $O(\ep)$ is uniform in $x\in \Omega$ with
$\mbox{\rm dist}(x,\cup_{i=1}^n\partial\Omega_i) \ge c$ for some constant $c>0$.
\end{theorem}

Note that
$$
\frac{1}{n}\sum_{i=1}^n \Phi_\Omega(x,x(s_i)) \rightarrow
v_\Omega(x):= \frac{1}{|\partial \Omega|} \int_{\partial \Omega}
\Phi_\Omega(x,z)\, d\sigma(z) \quad \mbox{as } n=
O(\frac{1}{\ep})\rightarrow + \infty,
$$
and the corrector $v_\Omega$ satisfies the equation
$$
\Delta v_\Omega = -1 \quad \mbox{in } \Omega,
$$
with, as already known in the homogenization theory
\cite{gadylshin},  the effective Neumann boundary condition
$$
\frac{\partial v_\Omega}{\partial \nu} = -
\frac{|\Omega|}{|\partial \Omega|} \quad \mbox{on } \partial
\Omega.
$$

\section{Asymptotics of the eigenvalues}

Let $0= \lambda_0^{(1)} < \lambda_0^{(2)} \le  \lambda_0^{(3)}
\leq \ldots$ be the eigenvalues of $-\Delta$ on $\Omega$ with
Neumann conditions on $\partial \Omega$. Let $u_0^{(j)}$ denote
the normalized eigenfunction associated with $\lambda_0^{(j)} $;
that is, it satisfies $\int_\Omega |u_0^{(j)}|^2 =1$. Let $\omega
\notin \{\sqrt{\lambda_0^{(j)}}\}_{j\geq 1}$. Introduce
$N^\omega(x,z)$ as the Neumann function for $\Delta + \omega^2$ in
$\Omega$ corresponding to a Dirac mass at $z$. That is, $N^\omega$
is the unique solution to
 \begin{equation} \label{defNeum}
 \left\{
 \begin{array}{ll}
 \ds(\Delta_{x} + \omega^2)
 N^\omega(x,z)=-\delta_{z}\quad &\mbox{in }\Omega,\\
 \nm\ds\frac{\partial N^\omega}{\partial \nu} =0\quad &\mbox{on
 }\partial \Omega.
 \end{array}
 \right.
\end{equation}
We recall two useful facts on the Neumann function \cite[(2.28) \& Lemma 5.1]{book}.

\begin{lemma} \label{neu-spec}
The following spectral decomposition holds pointwise:
\begin{equation} \label{spectrN}
N^\omega(x,z)=\sum_{j=1}^{+\infty}\frac{u^{(j)}_0(x)
{{u^{(j)}_0}(z)}}{\lambda^{(j)}_0 -\omega^2}, \quad x \neq z \in
\Omega.
\end{equation}
Moreover, if $\lambda^{(j_0)}_0$ is simple and $V$ is a complex neighborhood of $\sqrt{\lambda^{(j_0)}_0}$ which contains no $\sqrt{\lambda^{(j)}_0}$ other than $\sqrt{\lambda^{(j_0)}_0}$, then for $\omega \in V$ the Neumann function $N^\omega$ on the boundary $\partial\Omega$, which we denote by $N_{\partial \Omega}^\omega$,
has the form
\begin{equation} \label{singN} N_{\partial \Omega}^\omega (x,z) = - \frac{1}{\pi} \ln
|x-z| + \frac{u^{(j)}_0(x) {{u^{(j_0)}_0}(z)}}{\lambda^{(j_0)}_0
-\omega^2} + R_{\partial \Omega}^\omega(x,z) \quad \mbox{for } x
\in \Omega, z \in
\partial \Omega,
\end{equation} where $R_{\partial \Omega}^\omega(\cdot,z)$ belongs to $H^{3/2}(\Omega)$ for
any $z\in \partial \Omega$ and $\omega \mapsto R_{\partial
\Omega}^\omega$ is holomorphic in a complex neighborhood of
$\lambda_0^{(j_0)}$.
\end{lemma}

Suppose that $\lambda^{(j_0)}_0$ is  simple. Then there exists a
simple eigenvalue $\lambda^{(j_0)}_\ep$ near $\lambda^{(j_0)}_0$
associated to the normalized eigenfunction $u^{(j_0)}_{\ep}$
($\int_\Omega |u_\ep^{(j_0)}|^2 =1$) satisfying the eigenvalue
problem
\begin{equation}\label{PoissonEigen}
\begin{cases}
- \Delta u^{(j_0)}_{\ep}  = \lambda^{(j_0)}_\ep u^{(j_0)}_\ep & \quad  \mbox{in}~\Omega,\\
\ds u^{(j_0)}_{\ep}=0 & \quad \mbox{on}~{\partial \Omega_a},\\
\ds \frac{\partial u^{(j_0)}_{\ep}}{\partial \nu}=0 & \quad
\mbox{on}~{\partial \Omega_r}.
\end{cases}
\end{equation}

Our aim in this section is to derive a high-order asymptotic
expansion of $\lambda^{(j_0)}_\ep-\lambda^{(j_0)}_0$ as
$\ep\rightarrow 0$.

Set $\varphi_\ep = \partial u_\ep^{(j_0)}/\partial \nu$ on
$\partial \Omega_a$. By Green's formula, one can easily see that
$$
u_\ep(x) = (\lambda_\ep^{(j_0)} - \omega^2) \int_\Omega
u^{(j_0)}_\ep(z) N^\omega(x,z)\, dz + \int_{\partial \Omega_a}
\varphi_\ep(z) N^\omega(x,z)\, d\sigma(z),
$$
which shows that
$$
\int_{\partial \Omega_a} \varphi_\ep(z)
N^{\sqrt{\lambda_\ep^{(j_0)}}}(x,z)\, d\sigma(z) = 0 \quad
\mbox{on } \partial \Omega_a.
$$
 Therefore, the eigenvalue problem
\eqnref{PoissonEigen} reduces to the study of the characteristic
values of the operator-valued function $\lambda \mapsto
\mathcal{A}_\ep(\lambda)$ given by
$$\mathcal{A}_\ep(\lambda)[\varphi](x(t)) = \int_{-\ep}^\ep \varphi(s) N^{\sqrt{\lambda}}(x(t),x(s))\,
ds,
$$
where $x(t)$ is defined by \eqnref{paramet}.

Let $V$ be a neighborhood of $\lambda_0^{(j_0)}$ in the complex
plane such that $\lambda_0^{(j_0)}$ is the only eigenvalue of
$-\Delta$ in $\Omega$ with Neumann boundary condition. The following pole-pencil decomposition of
$\mathcal{A}_\ep(\lambda): \mathbf{X}_\ep \rightarrow
\mathbf{Y}_\ep$ holds for any $\lambda \in V$:
$$
\mathcal{A}_\ep(\lambda) = - \frac{1}{\pi} \mathcal{L}_\ep +
\frac{\mathcal{K}_\ep}{\lambda_0^{(j_0)} - \lambda} +
\mathcal{R}_\ep(\lambda),
$$
where
$$
{\mathcal L}_\ep[\varphi](x(t)) = \ds \int_{- \epsilon}^{\epsilon}
\ln|s-t|\,\varphi(s)\,ds,
$$
${{\mathcal K}}_\ep$ is the one-dimensional operator given by
\[
\ds {{\mathcal K}}_\ep[\varphi](x(t)) = \bigg(\int_{-
\epsilon}^{\epsilon} u_0^{j_0}(x(s)) \varphi(x(s))\, ds \bigg) \,
u_0^{j_0}(x(t)),
\]
and
\[
{{\mathcal R}}_\ep(\lambda)[\varphi](x(t)) = \ds \int_{-
\epsilon}^{\epsilon} R^{\sqrt{\lambda}}_{\partial \Omega}(x(t),
x(s)) \,\varphi(s)\,ds.
\]
An application of the generalized Rouch\'e theorem \cite{book}
shows that the operator-valued function $\mathcal{A}_\ep(\lambda)$
has exactly one characteristic value in $V$. Theorem 5.7 in
\cite{book} gives its asymptotic expansion for $\ep$ small enough.
The following theorem holds.
\begin{theorem} We have
 \beq \label{asymp-eigen}
 \begin{array}{lll} \ds
\lambda_\ep^{(j_0)} - \lambda_0^{(j_0)} &\approx&\ds - \pi
\int_{-\ep}^\ep \mathcal{L}_\ep^{-1}[u_{0}^{(j_0)}](x(t))
u_{0}^{(j_0)}(x(t))\, dt \\ && + \ds \pi^2 \int_{-\ep}^\ep
\mathcal{L}_\ep^{-1} {{\mathcal R}}_\ep(\lambda_0^{(j_0)})
\mathcal{L}_\ep^{-1} [u_{0}^{(j_0)}](x(t)) u_{0}^{(j_0)}(x(t))\,
dt. \end{array}
 \eeq
In particular, the leading-order term in the expansion of
$\lambda_\ep^{(j_0)} - \lambda_0^{(j_0)}$ is given by \beq
\lambda_\ep^{(j_0)} - \lambda_0^{(j_0)} \approx - \frac{\pi}{\ln
\ep} |u_0^{(j_0)}(x^*)|^2, \eeq where $x^*$ is the center of
$\partial \Omega_a$. \end{theorem}

Note that if $\lambda_0^{(j_0)}$ is a multiple eigenvalue of
$-\Delta$ in $\Omega$ with Neumann boundary condition, then in
exactly the same way as Theorem 3.24 in \cite{book}, we can
address the splitting problem in the evaluation of eigenvalues of
\eqnref{PoissonEigen}.

\section{Narrow escape problem in a field of force}

We now deal with the case when the force field $F=\nabla \phi$ is not zero. As before,
we consider the problem of estimating the escape time and that of deriving asymptotics
of the perturbed eigenvalue.

\subsection{Asymptotics of the solution}
Let the Neumann function  $N^F(x,z)$ be defined by
\begin{equation}\label{NeumannF}
\begin{cases}
&\Delta_z N^F(x,z) - \nabla_z \cdot (F(z) N^F(x,z)) =-\delta_x, \ \ x,z\in \Omega, \\
& \displaystyle\frac{\partial N^F}{\partial \nu_z} - F\cdot \nu_z
N^F \bigg|_{z \in
\partial \Omega} = - \frac{1}{|\partial \Omega|}.
\end{cases}
\end{equation}

The (boundary) Neumann function, denoted by $N_{\p\Om}^F(x,z)$, has the form
\begin{equation}
\label{Neumann smoothF}N_{\partial \Omega}^F(x,z)=-\frac{1}{\pi}
\ln|x-z| + R^F_{\partial\Omega}(x,z), \ \ x\in \Omega,\
z\in\partial\Omega,
\end{equation}
where $R^F_{\partial\Omega}(\cdot,z)$ belongs to
$H^{{3}/{2}}(\Omega)$ for any $z\in\partial\Omega$.

The solution $u_\ep$ to \eqnref{Poisson} admits the following
integral representation:
$$u_\ep(x)= \int_\Omega N^F(x,z) \, dz + \int_{\partial \Omega_a}
\varphi_\ep(z) N_{\p\Om}^F(x,z)\,d\sigma(z) + \frac{1}{|\partial \Omega|}
\int_{\p\Omega_r} u_\ep(z)\, d\sigma(z), \quad x \in \Omega,$$ where
$\varphi_\ep =
\partial u_\ep/\partial \nu$ on $\partial \Omega_a$. Therefore,
the density function $\varphi_\ep$ satisfies the integral equation
$$
- \frac{1}{\pi} \int_{\partial \Omega_a} \ln |x-z|
\varphi_\ep(z)\, d\sigma(z) + \int_{\partial \Omega_a} R_{\partial
\Omega}^F(x,z) \varphi_\ep(z)\, d\sigma(z) = - \int_\Omega
N^F(x,z)\, dz - C_\ep
$$
for $x \in \partial \Omega_a$, where
 $$
 C_\ep = \frac{1}{|\partial\Omega|} \int_{\p\Omega_r} u_\ep(z)\, d\sigma(z).
 $$
Moreover, since $F=\nabla \phi$, we have
$$\nabla \cdot (e^{\phi}\nabla  u_\ep) = - e^{\phi} \quad \mbox{in }
\Omega,$$
and hence $\varphi_\ep$ satisfies the compatibility condition
$$
\int_{\partial \Omega} e^{\phi(z)} \varphi_\ep(z)\, d\sigma(z) = -
\int_{\Omega} e^{\phi(z)} \, dz.
$$

In exactly the same way as Theorem \ref{themF0}, we can prove that
$$\begin{array}{lll}
\varphi_\ep(x(t)) &=&\ds \bigg[ e^{- \phi(x^*)} \int_\Omega
e^{\phi(z)}\, dz\, \bigg(\frac{\ln\ep}{\pi}-R_{\partial
\Omega}^F(x^*,x^*)\bigg) \\ \nm && \qquad \ds + \Ce+ \int_\Omega
N^F(x^*, z)\, dz\bigg] \frac{1}{(\ln\frac{1}{2})\, \sqrt{1-t^2}}+
O(\ep), \end{array}
$$
 where $x(t)$ is defined by
\eqnref{paramet} and $x^*$ is the center of $\partial \Omega_a$.
Consequently, we have
$$
C_\ep = \bigg(\int_\Omega e^{\phi(z)- \phi(x^*)}\,
dz\, \bigg) \bigg(\frac{1}{\pi} \ln{\frac{2}{\ep}} + R_{\partial
\Omega}^F(x^*,x^*)\bigg) - \int_\Omega N^F(x^*, z)\, dz + O(\ep),
$$
which yields the following result.
\begin{theorem}
The following asymptotic formula
$$
u_\ep(x)=\ds \frac{1}{\pi} \bigg(\int_\Omega e^{\phi(z)- \phi(x^*)} \bigg)
dz\, \ln\frac{2}{\ep}+\Phi^F_\Omega(x,x^*) +O(\ep),
$$
holds uniformly for $x\in \Omega$ satisfying $\mbox{\rm dist}(x,
\partial \Omega_a) \ge c$ for some constant $c>0$. Here
$\Phi^F_\Omega$ is given by
 \beq
 \begin{array}{rl}
 \Phi^F_\Omega(x,x^*) &=\ds\int_\Omega N^F(x,z)dz -\int_\Omega N^F(x^*,z)dz\\
 \nm & \qquad \ds + \bigg(\int_\Omega e^{\phi(z)- \phi(x^*)}\, dz\,
 \bigg) \big( R^F_{\partial\Omega}(x^*,x^*) - N^F_{\partial\Omega}(x,x^*) \big).
 \end{array}
 \eeq
\end{theorem}

\subsection{Asymptotics of the eigenvalues}
Consider the eigenvalue problem
\begin{equation}\label{PoissonEigenF}
\begin{cases}
&\Delta u_{\ep,F} + F \cdot \nabla u_{\ep,F}  = - \lambda_{\ep,F} u_{\ep,F}, \qquad  \mbox{in}~\Omega\\
&u_{\ep,F} =0 \qquad \mbox{on}~{\partial \Omega_a},\\
&\ds \frac{\partial u_{\ep,F}}{\partial \nu}=0 \qquad
\mbox{on}~{\partial \Omega_r}.
\end{cases}
\end{equation}
Since $F= \nabla \phi$ for a smooth potential function $\phi$, we
get
$$
e^{-\phi} \nabla \cdot (e^{\phi} \nabla u) = \Delta u +  F\cdot
\nabla u,$$ and hence we can realize the Neumann Laplacian with a
drift $F$ in $\Omega$ as a non-negative, self-adjoint operator on
$L^2(\Omega, e^{-\phi} dx)$ via the closure of the Friedrichs
extension of the nonnegative quadratic form
 $$
 (u,v) \mapsto \int_\Omega \nabla u \cdot \nabla v e^{\phi}\, dx.
 $$
Let $0 = \lambda_{0,F}^{(1)} < \lambda_{0,F}^{(2)} \le
\lambda_{0,F}^{(3)} \leq \ldots$ be the eigenvalues of $-\Delta -
F\cdot \nabla$ on $\Omega$ with Neumann conditions on $\partial
\Omega$ and let $u_{0,F}^{(1)},~ u_{0,F}^{(2)},~u_{0,F}^{(3)}, \ldots$ be the corresponding eigenfunctions with
$\int_\Om |u_{0,F}^{(j)}|^2e^\phi=1$. Define for $\omega \notin
\{\sqrt{\lambda_0^{(j)}}\}_{j\geq 1}$ the Neumann function by
\begin{equation}\label{NeumannFeigen}
\begin{cases}
&\Delta_z N^F_\omega(x,z) - \nabla_z \cdot (F(z) N^F_\omega(x,z))
+ \omega^2 N^F_\omega(x,z) =-\delta_x, \ \ x,z\in \Omega\\
& \displaystyle\frac{\partial N_\omega^F}{\partial \nu_z} - F\cdot
\nu_z N_\omega^F \bigg|_{z \in
\partial \Omega} =0.
\end{cases}
\end{equation}
Set

$$\varphi_{\ep,F}:= \frac{\partial u_{\ep,F}}{\partial \nu}
\quad \mbox{on } \partial \Omega_a. $$

From
$$
\int_{\partial \Omega_a} \varphi_{\ep,F}(z)
N^F_{\sqrt{\lambda_{\ep,F}}}(x,z)\, d\sigma(z) = 0 \quad \mbox{on
}
\partial \Omega_a,
$$
it follows that the eigenvalue problem \eqnref{PoissonEigenF}
reduces to the study of the characteristic values of the
operator-valued function $\lambda \mapsto
\mathcal{A}_{\ep,F}(\lambda)$ given by
$$\mathcal{A}_{\ep,F}(\lambda)[\varphi](x(t)) = \int_{-\ep}^\ep \varphi(s) N^F_{\sqrt{\lambda}}(x(t),x(s))\,
ds.
$$
Following the same arguments as for the case without drift, we can
prove the following result.

\begin{theorem} If $\lambda_{0,F}^{(j_0)}$ is simple then there exists a
simple eigenvalue $\lambda_{\ep,F}^{(j_0)}$ of
\eqnref{PoissonEigenF} converging to $\lambda_{0,F}^{(j_0)}$ as
$\ep$ goes to zero and the following formula holds: \beq
\lambda_{\ep,F}^{(j_0)} - \lambda_{0,F}^{(j_0)} \approx -
\frac{\pi}{\ln \ep} |u_{0,F}^{(j_0)}(x^*)|^2 e^{\phi(x^*)}. \eeq
\end{theorem}

\section{Conclusion}

In this paper, we have provided mathematically rigorous
derivations of the first- and second-order terms in the asymptotic
expansion of the solutions of diffusion equations in the presence
of a single or many small targets. We have shown the nonlinear
interaction between many targets. We have also studied the problem
of eigenvalue changes due to the small targets using the approach
developed in \cite{book}. The generalization of the results of
this work to nonsmooth domains containing corners and cusps will
be the subject of a forthcoming paper.

\section*{Acknowledgements} This work was supported  by R\'egion
\^{I}le-de-France and the grant NRF 20090085987.

\end{document}